\documentclass[onecolumn]{IEEEtran}

\usepackage{cite}
\usepackage[utf8]{inputenc}   % allow utf-8 input
\usepackage[T1]{fontenc}      % use 8-bit T1 fonts
\usepackage[hyphens]{url}     % simple URL typesetting
\usepackage{hyperref}         % hyperlinks
\usepackage{booktabs}         % professional-quality tables
\usepackage{amsmath}          % AMS math package
\usepackage{amsthm}           % AMS theorem package
\usepackage{amssymb}          % blackboard math symbols
\usepackage{nicefrac}         % compact symbols for 1/2, etc.
\usepackage{microtype}        % microtypography
\usepackage{enumerate}        % enumerate package for assumptions listing
\usepackage{algorithm}        % not supported by IEEEtran
\usepackage{algpseudocode}    % package for writing algorithms
\usepackage{pgfplots}         % LaTeX-native plots

\newtheorem{corollary}{Corollary}
\newtheorem{lemma}{Lemma}

\newtheorem{theorem}{Theorem}
\newtheorem{definition}{Definition}
\newtheorem{remark}{Remark}

\DeclareMathOperator*{\argmin}{arg\,min}

\newcommand{\dotprod}[2]{\ensuremath{\left\langle#1, #2\right\rangle}}

\newcommand{\norm}[1]{\ensuremath{\left\Vert#1\right\Vert}}

\definecolor{kth-dark-blue}{RGB/cmyk}{25,84,166/0.849,0.494,0,0.349}
\definecolor{kth-light-blue}{RGB/cmyk}{36,160,216/0.833,0.259,0,0.153}
\definecolor{kth-dark-red}{RGB/cmyk}{157,16,45/0,0.898,0.713,0.384}
\definecolor{kth-light-red}{RGB/cmyk}{228,54,62/0,0.763,0.728,0.106}
\definecolor{kth-dark-green}{RGB/cmyk}{98,146,46/0.329,0,0.685,0.427}
\definecolor{kth-light-green}{RGB/cmyk}{176,201,43/0.124,0,0.786,0.212}
\definecolor{kth-yellow}{RGB/cmyk}{250,185,25/0,0.26,0.9,0.0196}
\definecolor{kth-pink}{RGB/cmyk}{216,84,151/10,0.611,0.301,0.153}
\definecolor{kth-dark-grey}{RGB/cmyk}{101,101,108/0.0648,0.0648,0,0.576}
\definecolor{kth-grey}{RGB/cmyk}{189,188,188/0,0.00529,0.00529,0.259}
\definecolor{kth-light-grey}{RGB/cmyk}{227,229,227/0.00873,0,0.00873,0.102}

\pgfplotscreateplotcyclelist{kth-cycle}{%
  {thick,solid,kth-dark-blue},
  {thick,solid,kth-dark-red},
  {thick,solid,kth-dark-green},
  {thick,dashed,kth-light-blue},
  {thick,dashed,kth-light-red},
  {thick,dashed,kth-light-green},
  {thick,solid,kth-yellow},
  {thick,solid,kth-pink},
  {thick,solid,kth-dark-grey},
  {thick,solid,kth-grey},
  {thick,solid,kth-light-grey},
  {thick,dashed,kth-yellow},
  {thick,dashed,kth-pink},
  {thick,dashed,kth-dark-grey},
  {thick,dashed,kth-grey},
  {thick,dashed,kth-light-grey},
}
\pgfplotscreateplotcyclelist{kth-colors}{%
  {kth-dark-blue,fill=kth-dark-blue},
  {kth-dark-red,fill=kth-dark-red},
  {kth-dark-green,fill=kth-dark-green},
  {kth-light-blue,fill=kth-light-blue},
  {kth-light-red,fill=kth-light-red},
  {kth-light-green,fill=kth-light-green},
  {kth-yellow,fill=kth-yellow},
  {kth-pink,fill=kth-pink},
  {kth-dark-grey,fill=kth-dark-grey},
  {kth-grey,fill=kth-grey},
  {kth-light-grey,fill=kth-light-grey},
}

\begin{document}

\title{Analysis and Implementation of an Asynchronous Optimization Algorithm for
the Parameter Server}

\author{Arda~Aytekin~\IEEEmembership{Student~Member,~IEEE,}
Hamid~Reza~Feyzmahdavian,~\IEEEmembership{Student~Member,~IEEE,}
and~Mikael~Johansson~\IEEEmembership{Member,~IEEE}% <-this % stops a space
\thanks{A.~Aytekin, H.~R.~Feyzmahdavian and M.~Johansson are with the Department
of Automatic Control, School of Electrical Engineering and ACCESS Linnaeus
Center, KTH Royal Institute of Technology, SE-100 44 Stockholm, Sweden.
Emails: \texttt{\{aytekin, hamidrez, mikaelj\}@kth.se}}}

\maketitle

\begin{abstract}
  This paper presents an asynchronous incremental aggregated gradient algorithm
  and its implementation in a parameter server framework for solving regularized
  optimization problems. The algorithm can handle both general convex (possibly
  non-smooth) regularizers and general convex constraints. When the empirical
  data loss is strongly convex, we establish linear convergence rate, give
  explicit expressions for step-size choices that guarantee convergence to the
  optimum, and bound the associated convergence factors. The expressions have an
  explicit dependence on the degree of asynchrony and recover classical results
  under synchronous operation. Simulations and implementations on commercial
  compute clouds validate our findings.
\end{abstract}


\section{Introduction}

\IEEEPARstart{M}{achine} learning and optimization theory have enjoyed a 
%very
fruitful symbiosis over the last decade. On the one hand, since many machine learning
tasks can be posed as optimization problems, 
%so 
advances in 
%theory and
%algorithms for 
large-scale optimization (e.g.~\cite{Nes:12}) have had an
immediate and profound impact on machine learning research. On the other hand,
the challenges of dealing with huge data sets, often spread over multiple sites,
have inspired the machine learning community to develop novel 
%convex
optimization algorithms~\cite{DGS+:12}, improve the theory for asynchronous
computations~\cite{RRW+:11}, and introduce new programming models for parallel
and distributed optimization~\cite{LZY+:11}.

In this paper, we consider machine learning in the parameter server framework~\cite{LZY+:11}.
This is a master-worker architecture, where a central server maintains the
current parameter iterates and queries worker nodes for gradients of the loss
evaluated on their data. In this setting, we focus on problems on the form
\begin{align*}
  \begin{aligned}
    & \underset{x}{\text{minimize}} & &
    \sum_{n=1}^{N} f_{n}(x) + h(x) \\
    & \text{subject to}             & &
    x \in \mathbb{R}^{d} \,.
  \end{aligned}
\end{align*}
Here, the first part of the objective function typically models the empirical
data loss and the second term is a regularizer (for example, an $\ell_1$ penalty
to promote sparsity of the solution). Regularized optimization problems arise in
many applications in machine learning, signal processing, and statistical
estimation. Examples include Tikhonov and elastic net regularization, Lasso,
sparse logistic regression, and support vector machines.

In the parameter server framework, Li \emph{et al.}~\cite{LZY+:11} analyzed a
parallel and asynchronous proximal gradient method for non-convex problems and
established conditions for convergence to a critical point. Agarwal and
Duchi~\cite{AgD:11}, and more recently Feyzmahdavian \emph{et al.}~\cite{FAJ:16},
developed parallel mini-batch optimization algorithms based on asynchronous
incremental gradient methods. When the loss functions are strongly convex, which
is often the case, it has recently been observed that incremental aggregated
methods outperform incremental gradient descent and are, in addition, able to
converge to the true optimum even with a constant step-size. Gurbuzbalaban
\emph{et al.}~\cite{2015-Gurbuzbalaban} established linear convergence for an
incremental aggregated gradient method suitable for implementation in the
parameter server framework. However, the analysis does not allow for any
regularization term, nor any additional convex constraints. 

This paper presents an asynchronous proximal incremental aggregated gradient
algorithm and its implementation in the parameter server framework. Our
algorithm can handle both general convex regularizers and convex constraints.
We establish linear convergence when the empirical data loss is strongly convex,
give explicit expressions for step-size choices that guarantee convergence to
the global optimum and bound the associated convergence factors. These
expressions have an explicit dependence on the degree of asynchrony and recover
classical results under synchronous operation. We believe that this is a
practically and theoretically important addition to existing optimization
algorithms for the parameter server architecture.

\subsection{Prior work}

Incremental gradient methods for smooth optimization problems have a long
tradition, most notably in the training of neural networks via back-propagation.
In contrast to gradient methods, which compute the full gradient of the loss
function before updating the iterate, incremental gradient methods evaluate the
gradients of a single, or possibly a few, component functions in each iteration.
Incremental gradient methods can be computationally more efficient than
traditional gradient methods since each step is  cheaper but makes a comparable
progress on average. However, for global convergence, the step-size needs to
diminish to zero, which can lead to slow convergence~\cite{Bertsekas:11}. If a
constant step-size is used, only convergence to an approximate solution can be
guaranteed in general~\cite{Solodov:98}.

Recently, Blatt, Hero, and Gauchman~\cite{2007-Blatt} proposed a method,
the \textit{incremental aggregated gradient} (IAG), that also computes
the gradient of a single component function at each iteration. But rather
than updating the iterate based on this information, it uses the sum
of the most recently evaluated gradients of all component functions.
Compared to the basic incremental gradient methods, IAG has the advantage that
global convergence can be achieved using a constant step-size when each
component function is convex quadratic. Later, Gurbuzbalaban,
Ozdaglar, and Parillo~\cite{2015-Gurbuzbalaban} proved linear convergence for
IAG in a more general setting when component functions are strongly convex.
In a more recent work, Vanli, Gurbuzbalaban and Ozdaglar~\cite{2016-Gurbuzbalaban}
analyzed the global convergence rate of proximal incremental aggregated gradient
methods, where they can provide the linear convergence rate only after
sufficiently many iterations. Our result differs from theirs in that we provide
the linear convergence rate of the algorithm without any constraints on the
iteration count and we extend the result to the general distance functions.

There has been some recent work on the stochastic version of the IAG method
(called stochastic average gradient, or SAG) where we sample the component
function to update instead of using a cyclic
order~\cite{2016-Schmidt,2014-Defazio,Mairal:15}. Unlike the IAG method where
the linear convergence rate depends on the number of passes through the data,
the SAG method achieves a linear convergence rate that depends on the number of
iterations. Further, when the number of training examples is sufficiently large, 
the SAG method allows the use of a very large step-size, which leads to improved
theoretical and empirical performance.

\section{Notation}

We let $\mathbb{N}$ and $\mathbb{N}_{0}$ denote the set of natural numbers and
the set of natural numbers including zero, respectively. The inner product of
two vectors $x, y \in \mathbb{R}^{d}$ is denoted by $\dotprod{x}{y}$. We assume
that $\mathbb{R}^{d}$ is endowed with a norm $\norm{\cdot}$ and use
$\norm{\cdot}_{*}$ to represent the corresponding dual norm, defined by
\begin{align*}
  \norm{y}_{*} = \underset{\norm{x} \leq 1}{\sup} \dotprod{x}{y} \,.
\end{align*}

\section{Problem definition}

We consider optimization problems on the form
\begin{align}
\label{eqn:problem}
  \begin{aligned}
    & \underset{x}{\text{minimize}} & &
      \sum_{n=1}^{N} f_{n}(x) + h(x) \\
    & \text{subject to}             & &
    x \in \mathbb{R}^{d} \,,
  \end{aligned}
\end{align}
where $x$ is the decision variable, $f_{n}(x)$ is convex and differentiable for
each $n \in \mathcal{N} := \lbrace 1, \ldots, N \rbrace$ and $h(x)$ is a
proper convex function that may be non-smooth and extended real-valued. The role
of the regularization term $h(x)$ is to favor solutions with certain
preferred structure. For example, $h(x) = \lambda_{1} \norm{x}_{1}$ with
$\lambda_{1} > 0$ is often used to promote sparsity in solutions, and
\begin{align*}
  h(x) = I_{\mathcal{X}}(x) := \begin{cases}
    0         & \text{if } x \in \mathcal{X} \subseteq \mathbb{R}^{d} \,, \\
    +\infty   & \text{otherwise}
  \end{cases}
\end{align*}
is used to force the possible solutions to lie in the closed convex set
$\mathcal{X}$.

In order to solve~\eqref{eqn:problem}, we are going to use the proximal
incremental aggregated gradient method. In this method, at iteration $k
\in \mathbb{N}$, the gradients of all component functions $f_{n}(x)$, possibly
evaluated at stale information $x_{k-\tau_{k}^{n}}$, are aggregated 
\begin{align*}
  g_{k} = \sum_{n=1}^{N} \nabla f_{n} \left( x_{k-\tau_{k}^{n}} \right) \,.
\end{align*}
Then, a proximal step is taken based on the current vector $x_{k}$, the
aggregated gradient $g_{k}$, and the non-smooth term $h(x)$,
\begin{align}
\label{eqn:prox-euc}
  x_{k+1} = \argmin_{x} \bigg\lbrace \dotprod{g_{k}}{x-x_{k}} +
    \frac{1}{2\alpha} \norm{x - x_{k}}^{2} + h(x) \bigg\rbrace.
\end{align}
The algorithm has a natural implementation in the parameter server framework.
The master node maintains the iterate $x$ and performs the proximal steps.
Whenever a worker node reports new gradients, the master updates the iterate
and informs the worker about the new iterate. Pseudo code for a basic parameter
server implementation is given in Algorithms~\ref{alg:master}~and~\ref{alg:worker}.

\begin{algorithm}
  \caption{Master procedure}\label{alg:master}
  \begin{algorithmic}[1]
    \State \textbf{Data:} $g_{w}$ for each worker $w \in \mathcal{W} :=
      \lbrace 1, 2, \ldots, W \rbrace$

    \State \textbf{Input:} $\alpha(L, \mu, \bar{\tau})$, $K > 0$ and $h(x)$

    \State \textbf{Output:} $x_{K}$
    \State \textbf{Initialize:} $k = 0$

    \Statex

    \While{$k < K$}
      \State Wait until a set $\mathcal{R}$ of workers return their gradients

      \ForAll{$w \in \mathcal{W}$}
        \If{$w \in \mathcal{R}$}
          \State Update $g_{w} \gets \sum_{n \in \mathcal{N}_{w}} \nabla
            f_{n}(x_{k-\tau_{k}^{w}})$
        \Else
          \State Keep old $g_{w}$
        \EndIf
      \EndFor

      \State Aggregate the incremental gradients $g_{k} = \sum_{w \in
        \mathcal{W}} g_{w}$
      \State Solve~\eqref{eqn:prox-euc} with $g_{k}$

      \ForAll{$w \in \mathcal{R}$}
        \State Send $x_{k+1}$ to worker $w$
      \EndFor

      \State Increment $k$
    \EndWhile

    \State Signal \texttt{EXIT}
    \State Return $x_{K}$
  \end{algorithmic}
\end{algorithm}

\begin{algorithm}
  \caption{Procedure for each worker $w$}\label{alg:worker}
  \begin{algorithmic}[1]
    \State \textbf{Data:} $x$, and loss functions $
    	%\mathcal{F}_{w} = 
    \lbrace
      f_{w}(x) : w \in \mathcal{N}_{w} \rbrace$ with $\mathop{\bigcup}_{w \in
      \mathcal{W}} \mathcal{N}_{w} = \mathcal{N}$ and $\mathcal{N}_{w_{1}}
      \bigcap \mathcal{N}_{w_{2}} = \emptyset$ $\forall w_{1} \ne w_{2} \in
      \mathcal{W}$

    \Statex

    \Repeat
      \State Receive $x \gets x_{k+1}$ from master
      \State Calculate incremental gradient (IG) $\sum_{n \in \mathcal{N}_{w}} \nabla
        f_{n}(x)$
      \State Send IG to master with a delay of $\tau_{k}^{w}$
    \Until{\texttt{EXIT} received}
  \end{algorithmic}
\end{algorithm}

To establish convergence of the iterates to the global optimum, we impose the
following assumptions on Problem~\eqref{eqn:problem}:
\begin{enumerate}[{A}1)]
 \item\label{ass:first} The function $F(x) := \sum_{n=1}^{N} f_{n}(x)$ is $\mu$-strongly convex,
    \textit{i.e.},
    \begin{align}\label{eqn:strong-cvx}
      F(x) \geq F(y) + \dotprod{\nabla F(y)}{x-y} + \frac{\mu}{2}
        \norm{x-y}^{2},
    \end{align}
    holds for all $x, y \in \mathbb{R}^{d}$.
  \item Each $f_{n}$ is convex with $L_{n}$-continuous
    gradient, that is,
    \begin{align*}
      \norm{\nabla f_{n}(x) - \nabla f_{n}(y)}_{*} \leq L_{n} \norm{x-y} \quad
        \forall x, y \in \mathbb{R}^{d} \,.
    \end{align*}
Note that under this assumption, $\nabla F$ is also Lipschitz continuous with $L\leq \sum_{n=1}^{N} L_n$.
 \item $h(x)$ is sub-differentiable everywhere in its effective
    domain, that is, for all $x, y \in \lbrace z \in \mathbb{R}^{d} :  h(z) <
    +\infty \rbrace$,
    \begin{align}\label{eqn:subdiff}
      h(x) \geq h(y) + \dotprod{s(y)}{x-y} \quad \forall s(y) \in \partial h(y).
    \end{align}
  \item\label{ass:last} The time-varying delays $\tau_k^n$ are bounded,
    \textit{i.e.}, there is a non-negative integer $\bar{\tau}$ such that
    \begin{align*}
      \tau_{k}^{n} \in \lbrace 0, 1, \ldots, \bar{\tau} \rbrace,
    \end{align*}
hold for all $k \in \mathbb{N}_{0}$ and $ n \in \mathcal{N}$. 
\end{enumerate}

\section{Main result}

First, we provide a lemma which is key to proving our main result.

\begin{lemma}\label{lem:lemma}
  Assume that the non-negative sequences $\{ V_{k}\}$ and $\{w_{k}\}$ satisfy
  the following inequality:
  \begin{align}\label{eqn:sequence}
    V_{k+1} \leq a V_{k} - b w_{k} + c \sum_{j=k-k_{0}}^{k} w_{j} \,,
  \end{align}
  for some real numbers $a \in (0,1)$ and $b,c \geq 0$, and some integer $k_{0}
  \in \mathbb{N}_{0}$. Assume also that $w_{k} = 0$ for $k < 0$, and that the
  following holds:
  \begin{align*}
    \frac{c}{1-a} \frac{1-a^{k_{0}+1}}{a^{k_{0}}} \leq b \,.
  \end{align*}
  Then, $V_{k} \leq a^k V_0$ for all $k\geq 0$.
\end{lemma}

\begin{IEEEproof}
  To prove the linear convergence of the sequence, we divide both sides of~\eqref{eqn:sequence} by $a^{k+1}$ and take the sum:
  \begin{IEEEeqnarray}{rCl}
    \sum_{k=0}^{K} \frac{V_{k+1}}{a^{k+1}}  & \leq &
      \sum_{k=0}^{K} \frac{V_{k}}{a^{k}} - b\sum_{k=0}^{K} \frac{w_{k}}{a^{k+1}}
      + c \sum_{k=0}^{K} \frac{1}{a^{k+1}} \sum_{j=k-k_{0}}^{k} w_{j} \notag \\
                                            & = &
      \sum_{k=0}^{K} \frac{V_{k}}{a^{k}} - b\sum_{k=0}^{K} \frac{w_{k}}{a^{k+1}}
      \notag \\
                                            & &
      +\> \frac{c}{a} \left(w_{-k_{0}} + w_{-k_{0}+1} + \dots + w_{0}\right)
      \notag \\
                                            & &
      +\> \frac{c}{a^{2}} \left(w_{-k_{0}+1} + w_{-k_{0}+2} + \dots +
      w_{1}\right) + \dots \notag \\
                                            & &
      +\> \frac{c}{a^{K+1}} \left(w_{K-k_{0}} + w_{K-k_{0}+1}
      + \dots w_{K}\right) \notag \\
                                            & \leq &
      \left(c \left(1 + \frac{1}{a} +
      \dots + \frac{1}{a^{k_{0}}}\right) - b\right)
      \sum_{k=0}^{K} \frac{w_{k}}{a^{k+1}} \notag \\
                                            & &
      +\> \sum_{k=0}^{K} \frac{V_{k}}{a^{k}} \,, \label{eqn:seq-ineq}
  \end{IEEEeqnarray}
  where we have used the non-negativity of $w_{k}$ to
  obtain~\eqref{eqn:seq-ineq}.

  If the coefficient of the first sum of the right-hand side
  of~\eqref{eqn:seq-ineq} is non-positive, \textit{i.e.}, if
  \begin{align*}
    c + \frac{c}{a} + \dots + \frac{c}{a^{k_{0}}} =
      \frac{c}{1-a} \frac{1-a^{k_{0}+1}}{a^{k_{0}}} \leq b \,,
  \end{align*}
  then inequality~\eqref{eqn:seq-ineq} implies that
  \begin{align*}
    \frac{V_{K+1}}{a^{K+1}} + \frac{V_{K}}{a^{K}} + \dots +
      \frac{V_{1}}{a^{1}} \leq \frac{V_{K}}{a^{K}} +
      \frac{V_{K-1}}{a^{K-1}} \dots + \frac{V_{0}}{a^{0}}. %\,,
  \end{align*}
  Hence, $V_{K+1} \leq a^{K+1} V_{0}$ for any $K
  \geq 0$ and the desired result follows.
\end{IEEEproof}

We are now ready to state and prove our main result.

\begin{theorem}\label{thm:theorem}
  Assume that Problem~\eqref{eqn:problem} satisfies assumptions
  A\ref{ass:first}--A\ref{ass:last}, and that the step-size $\alpha$ satisfies:
  \begin{align*}%\label{eqn:step-size-1}
    \alpha \leq \frac{\left(  1+\frac{\mu}{L}\frac{1}{\bar{\tau}
      +1}\right)^\frac{1}{(\bar{\tau}+1)}
-1}{\mu} \,,
  \end{align*}
  where $L = \sum_{n=1}^{N} L_{n}$.  Then, the iterates generated by
  Algorithms~\ref{alg:master}~and~\ref{alg:worker} satisfy:
  \begin{align*}
    \norm{x_{k} - x^{\star}}^{2} \leq
      {\left(\frac{1}{\mu\alpha + 1}\right)}^{k}
      \norm{x_{0} - x^{\star}}^{2} \,.
  \end{align*}
  for all $k\geq 0$.
\end{theorem}

\begin{IEEEproof}
  We start with analyzing each component function $f_{n}(x)$ to find upper
  bounds on the function values:
  \begin{IEEEeqnarray}{rCl}
    f_{n}(x_{k+1})  & \leq &
      f_{n}(x_{k-\tau_{k}^{n}}) + \dotprod{\nabla f_{n}(x_{k-\tau_{k}^{n}})}
      {x_{k+1} - x_{k-\tau_{k}^{n}}} \notag \\
                    & &
      +\> \frac{L_{n}}{2} \norm{x_{k+1} - x_{k-\tau_{k}^{n}}}^{2} \notag \\
                    & \leq &
      f_{n}(x) + \dotprod{\nabla f_{n}(x_{k-\tau_{k}^{n}})}{x_{k+1} - x}
        \notag \\
                    & &
      +\> \frac{L_{n}}{2} \norm{x_{k+1} - x_{k-\tau_{k}^{n}}}^{2}
        \quad \forall x \,, \label{eqn:first-cvx-bnd}
  \end{IEEEeqnarray}
  where the first and second inequalities use $L_{n}$-continuity and convexity
  of $f_{n}(x)$, respectively. Summing~\eqref{eqn:first-cvx-bnd} over all
  component functions, we obtain:
  \begin{IEEEeqnarray}{rCl}
    F(x_{k+1})  & \leq & F(x) + \dotprod{g_{k}}{x_{k+1}-x} \notag \\
                & &
    +\> \sum_{n=1}^{N} \frac{L_{n}}{2} \norm{x_{k+1} -
      x_{k-\tau_{k}^{n}}}^{2} \quad \forall x \,.
      \label{eqn:second-cvx-bnd}
  \end{IEEEeqnarray}

  Observe that optimality condition of~\eqref{eqn:prox-euc} implies:
  \begin{IEEEeqnarray}{rCl}\label{eqn:optimality-euc}
    \dotprod{g_{k}}{x_{k+1} - x} & \leq &
      \frac{1}{\alpha} \dotprod{x_{k+1} - x_{k}}{x - x_{k+1}} \notag \\
      & &
      +\> \dotprod{s(x_{k+1})}{x - x_{k+1}} \quad \forall x \in \mathcal{X} \,.
  \end{IEEEeqnarray}

  To find an upper bound on the second term of the right-hand side
  of~\eqref{eqn:second-cvx-bnd}, we use the three-point equality
  on~\eqref{eqn:optimality-euc} to obtain:
  \begin{IEEEeqnarray}{rCl}\label{eqn:optimality-euc-three}
    \dotprod{g_{k}}{x_{k+1}-x} & \leq & \frac{1}{2\alpha} \norm{x_{k}-x}^{2}
      - \frac{1}{2\alpha} \norm{x_{k+1}-x_{k}}^{2} \notag \\
        & &
      -\> \frac{1}{2\alpha} \norm{x_{k+1}-x}^{2} \notag \\
        & &
      +\> \dotprod{s(x_{k+1})}{x-x_{k+1}} \quad \forall x \in \mathcal{X} \,.
  \end{IEEEeqnarray}

  Plugging $y = x_{k+1}$ in~\eqref{eqn:subdiff}, and using~\eqref{eqn:subdiff}
  together with~\eqref{eqn:optimality-euc-three} in~\eqref{eqn:second-cvx-bnd},
  we obtain the following relation:
  \begin{IEEEeqnarray*}{rCl}
    F(x_{k+1}) & + & h(x_{k+1}) + \frac{1}{2\alpha} \norm{x_{k+1}-x}^{2}
      \leq F(x) + h(x) \\
      & + & \frac{1}{2\alpha} \norm{x_{k}-x}^{2} -
      \frac{1}{2\alpha} \norm{x_{k+1}-x_{k}}^{2} \\
      & + & \sum_{n=1}^{N} \frac{L_{n}}{2}
      \norm{x_{k+1}-x_{k-\tau_{k}^{n}}}^{2}
      \quad \forall x \in \mathcal{X} \,.
  \end{IEEEeqnarray*}

  Using the strong convexity property~\eqref{eqn:strong-cvx} on $F(x_{k+1}) +
  h(x_{k+1})$ above and choosing $x = x^{\star}$ gives:
  \begin{IEEEeqnarray}{rCl}\label{eqn:third-cvx-bnd}
    \langle\nabla F(x^{\star}) & + & s(x^{\star}), x_{k+1} - x^{\star}\rangle
    + \frac{\mu}{2} \norm{x_{k+1}-x^{\star}}^{2} \notag \\
    & & +\> \frac{1}{2\alpha} \norm{x_{k+1}-x^{\star}}^{2} \notag \\
    & \leq & \frac{1}{2\alpha} \norm{x_{k}-x^{\star}}^{2} -
      \frac{1}{2\alpha} \norm{x_{k+1}-x_{k}}^{2} \notag \\
    & & +\> \sum_{n=1}^{N} \frac{L_{n}}{2} \norm{x_{k+1}-
      x_{k-\tau_{k}^{n}}}^{2} \,.
  \end{IEEEeqnarray}
  Due to the optimality condition of~\eqref{eqn:problem}, there exists a
  subgradient $s(x^{\star})$ such that the first term on the left-hand
  side is non-negative. Using this particular subgradient, we drop the first
  term. The last term on the right-hand side of the inequality can be further
  upper-bounded using Jensen's inequality as follows:
  \begin{IEEEeqnarray*}{rCl}
    \sum_{n=1}^{N} \frac{L_{n}}{2} \norm{x_{k+1}-x_{k-\tau_{k}^{n}}}^{2}
      & = &
      \sum_{n=1}^{N} \frac{L_{n}}{2} \norm{\sum_{j=k-\tau_{k}^{n}}^{k} x_{j+1} -
        x_{j}}^{2} \\
      & \leq & \frac{L(\bar{\tau}+1)}{2} \sum_{j=k-\bar{\tau}}^{k}
        \norm{x_{j+1}-x_{j}}^{2} \,,
  \end{IEEEeqnarray*}
  where $L = \sum_{n=1}^{N} L_{n}$. As a result, rearranging the terms
  in~\eqref{eqn:third-cvx-bnd}, we obtain:
  \begin{IEEEeqnarray*}{rCl}
    \norm{x_{k+1}-x^{\star}}^{2} & \leq &
      \frac{1}{\mu\alpha+1} \norm{x_{k}-x^{\star}}^{2} \notag \\
      & & -\> \frac{1}{\mu\alpha+1} \norm{x_{k+1}-x_{k}}^{2} \notag \\
      & & +\> \frac{\alpha(\bar{\tau}+1)L}{\mu\alpha+1}
        \sum_{j=k-\bar{\tau}}^{k} \norm{x_{j+1}-x_{j}}^{2} \,.
  \end{IEEEeqnarray*}

  We note that $\norm{x_{j+1}-x_{j}}^{2} = 0$ for all $j < 0$. Using
  Lemma~\ref{lem:lemma} with $V_{k} = \norm{x_{k+1}-x^{\star}}^{2}$, $w_{k}
  = \norm{x_{k+1}-x_{k}}^{2}$, $a = b = \frac{1}{\mu\alpha+1}$, $c =
  \frac{\alpha(\bar{\tau}+1)L}{\mu\alpha+1}$ and $k_{0} = \bar{\tau}$ completes
  the proof.
\end{IEEEproof}

\begin{remark}
	For the special case of Algorithms~\ref{alg:master}~and~\ref{alg:worker} where
$\tau_{k}^{n} = 0$ for all $k, n$, Xiao and Zhang~\cite{Xiao:14} have shown that
the convergence rate of serial proximal gradient method with a constant
step-size $\alpha = \frac{1}{L}$ is
\begin{align*}
\mathcal{O} \left( 
\left( \frac{L - \mu_{F}}{L + \mu_{h}} \right)^{k} \right)
\end{align*}
where $\mu_{F}$ and $\mu_{h}$ are strong convexity parameters of $F(x)$ and
$h(x)$, respectively. It is clear that in the case that $\bar{\tau}=0$, the
guaranteed bound in Theorem 1 reduces to the one obtained in~\cite{Xiao:14}.
\end{remark}

\section{Proximal incremental aggregate descent with general distance functions}

The update rule of our algorithm can be easily extended to a non-Euclidean
setting, by replacing the Euclidean squared distance in~\eqref{eqn:prox-euc}
with a general Bregman distance function. We first define a Bregman distance
function, also referred to as a prox-function.
\begin{definition}
\textup{A function $\omega \colon \mathbb{R}^{d} \to \mathbb{R}$ is called a
distance generating function with modulus $\mu_{\omega} > 0$ if $\omega$ is
continuously differentiable and $\mu_{\omega}$-strongly convex with respect to
$\norm{\cdot}$. Every distance generating function introduces a corresponding
Bregman distance function given by
\begin{align*}
  D_{\omega}(x,x') := \omega(x') - \omega(x) - \dotprod{\nabla\omega(x)}{x'-x}
    \,.
\end{align*}%
}
\end{definition}
For example, if we choose $\omega(x) = \frac{1}{2}\norm{x}_{2}^{2}$, which is
$1$-strongly convex with respect to the $l_{2}$-norm, that would result in
$D_{\omega}(x,x') = \frac{1}{2}\norm{x'-x}_{2}^{2}$. Another common example of
distance generating functions is the entropy function
\begin{align*}
  \omega(x) = \sum_{i=1}^{d} x_{i} \log(x_{i}) \,,
\end{align*}
which is $1$-strongly convex with respect to the $l_1$-norm over the standard
simplex
\begin{align*}
  \Delta := \left\lbrace x \in \mathbb{R}^{d} : \sum_{i=1}^{d} x_{i} = 1, \;
    x \geq 0 \right\rbrace \,,
\end{align*}
and its associated Bregman distance function is
\begin{align*}
  D_{\omega}(x,x') = \sum_{i=1}^{d} x'_{i} \log\left(\frac{x'_{i}}{x_{i}}\right)
    \,.
\end{align*}
The main motivation to use a generalized distance generating function rather
than the usual Euclidean distance function is to design an optimization
algorithm that can take advantage of the geometry of the feasible set.

The associated convergence result now reads as follows.
\begin{corollary}\label{cor:corollary}
  Consider using the following proximal gradient method to
  solve~\eqref{eqn:problem}:
  \begin{IEEEeqnarray}{rCl}\label{eqn:prox-gen}
      x_{k+1} & = & \argmin_{x \in \mathcal{X}} \bigg\lbrace
        \dotprod{g_{k}}{x-x_{k}} + \frac{1}{\alpha} D_{\omega}(x, x_{k})
        + h(x) \bigg\rbrace \,, \notag \\
      g_{k}   & = & \sum_{n=1}^{N} \nabla f_{n}
        \left(x_{k-\tau_{k}^{n}}\right) \,.
  \end{IEEEeqnarray}
  Assume that $D_{\omega}(\cdot, \cdot)$ satisfies:
  \begin{align}\label{eqn:gen-dist-bnd}
    \frac{\mu_{\omega}}{2} \norm{x-y}^{2} \leq D_{\omega}(x,y) \leq
      \frac{L_{\omega}}{2} \norm{x-y}^{2} \,.
  \end{align}
  Assume also that the problem satisfies assumptions
  A\ref{ass:first}--A\ref{ass:last}, and that the step-size $\alpha$ satisfies:
  \begin{align*}%\label{eqn:step-size-2}
    \alpha \leq \frac{
    	L_{\omega}
    	\left(
    	1+\frac{\mu}{L}
	      \frac{1}{\bar{\tau}+1}\frac{\mu_{\omega}}{L_{\omega}} \right)^{\frac{1}{\bar{\tau}+1}}-1 }{\mu} \,,
  \end{align*}
  where $L = \sum_{n=1}^{N} L_{n}$. Then, the iterates generated by the method
  satisfy:
  \begin{align*}
    D_{\omega}(x^{\star},x_{k}) \leq
      {\left(\frac{L_{\omega}}{\mu\alpha + L_{\omega}}\right)}^{k}
      D_{\omega}(x^{\star},x_{0}) \,.
  \end{align*}
\end{corollary}

\begin{IEEEproof}
  The analysis is similar to that of Theorem~\ref{thm:theorem}. This time, the
  optimality condition of~\eqref{eqn:prox-gen} implies:
  \begin{IEEEeqnarray}{rCl}\label{eqn:optimality-gen}
    \dotprod{g_{k}}{x_{k+1} - x} & \leq &
      \frac{1}{\alpha} \dotprod{\nabla\omega(x_{k+1}) - \nabla\omega(x_{k})}{x
        - x_{k+1}} \notag \\
      & &
      +\> \dotprod{s(x_{k+1})}{x - x_{k+1}} \quad \forall x \in \mathcal{X} \,.
  \end{IEEEeqnarray}

  Using the following four-point equality
  \begin{IEEEeqnarray*}{rCl}
    D_{\omega}(a,d) & - & D_{\omega}(c,d) - D_{\omega}(a,b) + D_{\omega}(c,b) =
      \notag \\
    & & \dotprod{\nabla\omega(b) - \nabla\omega(d)}{a-c} \,,
  \end{IEEEeqnarray*}
  in~\eqref{eqn:optimality-gen} with $a = x$, $b = c = x_{k+1}$ and $d = x_{k}$,
  and following the steps of the proof of Theorem~\ref{thm:theorem}, we obtain:
  \begin{IEEEeqnarray*}{rCl}
    \frac{\mu}{2} \lVert x_{k+1} & - & x^{\star} \rVert^{2} +
      \frac{1}{\alpha} D_{\omega}(x^{\star},x_{k+1}) \\
    & \leq & \frac{1}{\alpha} D_{\omega}(x^{\star},x_{k}) -
      \frac{1}{\alpha} D_{\omega}(x_{k+1},x_{k}) \\
    & & +\>
      \frac{L(\bar{\tau}+1)}{2} \sum_{j=k-\bar{\tau}}^{k}
      \norm{x_{j+1} - x_{j}}^{2} \,.
  \end{IEEEeqnarray*}

  This time, using the upper and lower bounds of~\eqref{eqn:gen-dist-bnd} on the
  left and right hand-side of the above inequality, respectively, and
  rearranging the terms, we arrive at:
  \begin{IEEEeqnarray*}{rCl}
    D_{\omega}(x^{\star},x_{k+1}) & \leq &
      \frac{L_{\omega}}{\mu\alpha + L_{\omega}} D_{\omega}(x^{\star},x_{k}) \\
    & & -\> \frac{L_{\omega}}{\mu\alpha+L_{\omega}} D_{\omega}(x_{k+1},x_{k})\\
    & & +\> \frac{\alpha L(\bar{\tau}+1)}{\mu\alpha+L_{\omega}}
      \frac{L_{\omega}}{\mu_{\omega}} \sum_{j=k-\bar{\tau}}^{k}
      D_{\omega}(x_{k+1},x_{k})
  \end{IEEEeqnarray*}

  Applying Lemma~\ref{lem:lemma} with $V_{k} = D_{\omega}(x^{\star},x_{k+1})$,
  $w_{k} = D_{\omega}(x_{k+1},x_{k})$, $a = b = \frac{L_{\omega}}{\mu\alpha +
  L_{\omega}}$, $c = \frac{\alpha L (\bar{\tau}+1)}{\mu\alpha + L_{\omega}}
  \frac{L_{\omega}}{\mu_{\omega}}$ and $k_{0} = \bar{\tau}$ completes the proof.
\end{IEEEproof}

\section{Numerical examples}

In this section, we present numerical examples which verify our theoretical
bound in different settings. First, we simulate the implementation of
Algorithms~\ref{alg:master}~and~\ref{alg:worker} on a parameter server
architecture to solve a small, toy problem. Then, we %will 
implement the framework on Amazon EC2 and solve a binary classification problem on three
different real-world datasets.

\subsection{Toy problem}

To verify our theoretical bounds provided in Theorem~\ref{thm:theorem} and
Corollary~\ref{cor:corollary}, we consider solving~\eqref{eqn:problem} with
\begin{align*}
  f_{n}(x)    & = \begin{cases}
    {(x_{n}-c)}^{2} + \frac{1}{2} {(x_{n+1}+c)}^{2}
      \,, \quad n = 1 \,, \\
    \frac{1}{2} {(x_{n-1}+c)}^{2} + \frac{1}{2} {(x_{n}-c)}^{2}
      \,, \quad n = N \,, \\
    \frac{1}{2} {(x_{n-1}+c)}^{2} + \frac{1}{2} {(x_{n}-c)}^{2} +
      \frac{1}{2} {(x_{n+1}+c)}^{2} \,,
  \end{cases} \\
  h(x)        & = \lambda_{1} \norm{x}_{1} + I_{\mathcal{X}}(x) \,, \\
  \mathcal{X} & = \left\lbrace x \geq 0 \right\rbrace \,,
\end{align*}
for some $c \geq 0$. We use $D_{\omega}(x,x_{k}) = \frac{1}{2}
\norm{x-x_{k}}_{p}^{2}$ in the proximal step~\eqref{eqn:prox-gen} and consider
both $p=1.5$ and $p=2$.

It can be verified that $\nabla F(x)$ is $(N+1)$-continuous and $F(x)$ is
2-strongly convex, both with respect to $\norm{\cdot}_{2}$, and that the
optimizer for the problem is $x^{\star} = \frac{\max(0, c-\lambda_{1})}{3}
e_{1}$, where $e_{n}$ denotes the $n^{\text{th}}$ basis vector. Moreover,
it can be shown that if $p\in (1, 2]$, then $\mu_{\omega} = 1$ and $L_{\omega}
= N^{2/p-1}$ satisfy~\eqref{eqn:gen-dist-bnd} with respect to $\norm{\cdot}_{2}$.

\begin{figure}[ht]
  \centering
  \begin{tikzpicture}
    \begin{semilogyaxis}[%
    xlabel=$k$,%
    ylabel=$\norm{x_{k}-x^{\star}}_{2}^{2}$,%
    grid=both,%
    scaled x ticks=false,%
    xmin=0,xmax=15000,%
    ymin=1e-9,ymax=4e3,%
    xtick={0,3000,6000,9000,12000,15000},%
    title=Iterate Convergence,%
    legend style={at={(0.5,-0.20)},anchor=north,legend columns=-1},%
    ]
    % p = 1.5 bound
    \addplot[kth-dark-blue,solid,thick]
      table[x index=0,y index=2,col sep=comma,header=true]
      {toy-results.csv};
    % p = 2.0 bound
    \addplot[kth-dark-red,solid,thick]
      table[x index=0,y index=4,col sep=comma,header=true]
      {toy-results.csv};
    % and their simulation results
    \addplot[kth-dark-blue,dashdotted,thick]
      table[x index=0,y index=1,col sep=comma,header=true]
      {toy-results.csv};
    \addplot[kth-dark-red,dashdotted,thick]
      table[x index=0,y index=3,col sep=comma,header=true]
      {toy-results.csv};
    % legend entries
    \legend{$p=1.5$,$p=2.0$}
    \end{semilogyaxis}
  \end{tikzpicture}
  \caption{Convergence of the iterates in toy problem. Solid lines represent
    our theoretical upper bound, whereas dash-dotted lines represent simulation
    results.}\label{fig:toy}
\end{figure}
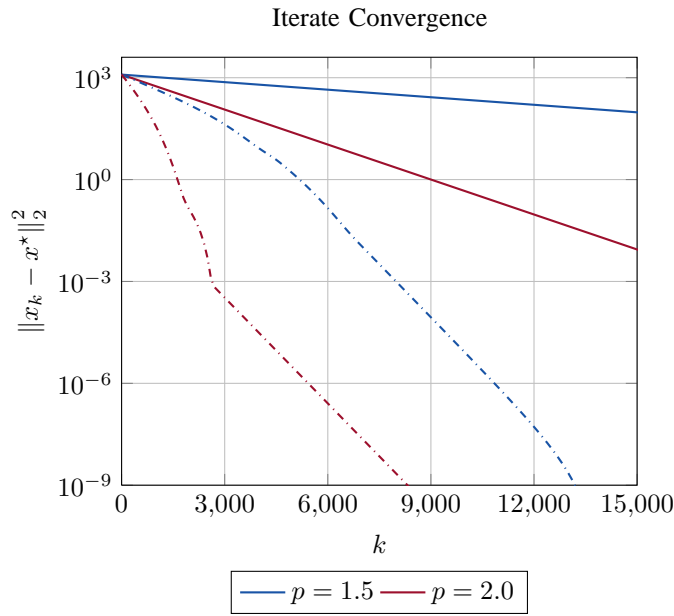

We select the problem parameters $N = 100$, $c = 3$ and $\lambda_{1} = 1$. We
simulate solving the problem with $W = 4$ workers, where at each iteration
$k$, a worker $w$ is selected uniformly at random to return their gradient, 
information evaluated on stale information $x_{k-\tau_{k}^{w}}$, to the master.
Here, at time $k$, $\tau_{k}^{w}$ is simply the number of iterations since the
last time worker $w$ was selected. Each worker holds $N/W=25$ component
functions, and we tune step-size based on the assumption that $\bar{\tau}=W$.
Figure~\ref{fig:toy} shows the results of a representative simulation. As can
be observed, the iterates converge to the optimizer and the theoretical bound
derived is valid.

\subsection{Binary classification on actual datasets}

Next, we consider solving a regularized, sparse binary classification
problem on three different datasets: \texttt{rcv1} (sparse)~\cite{2004-Lewis},
\texttt{url} (sparse)~\cite{URL:09} and \texttt{epsilon}
(dense)~\cite{Epsilon:09}. To this end, we implement the parameter server
framework in the Julia language, and instantiate it with Problem~\eqref{eqn:problem}:
\begin{align*}
  f_{n}(x)  & = \frac{1}{N} \left( \log(1+\exp(-b_{n}\dotprod{a_{n}}{x})) +
  \frac{1}{2} \lambda_{2} \norm{x}_{2}^{2} \right) \,, \\
  h(x)      & = \lambda_{1} \norm{x}_{1} \,,
\end{align*}
Here,  $a_{n} \in \mathbb{R}^{d}$ is the feature vector for  sample $n$, and
$b_{n} \in \lbrace -1, 1 \rbrace$ is the corresponding binary label. We pick
$\lambda_{1} = 10^{-5}$ and $\lambda_{2} = 10^{-4}$ for \texttt{rcv1} and
\texttt{epsilon} datasets, and $\lambda_{1} = 10^{-3}$ and $\lambda_{2} =
10^{-4}$ for \texttt{url}. \texttt{rcv1} is already normalized to have unit
norm in its samples; hence, we normalize \texttt{url} and \texttt{epsilon}
datasets to have comparable problem instances.

\texttt{rcv1} is a text categorization test collection from Reuters, having
$N = 804414$ documents and $d = 47236$ features (density: $0.16\%$) for
each document. We choose to classify \emph{sports}, \emph{disaster} and
\emph{government} related articles from the corpus. %, \textit{i.e.}, $b_{n} = 1$.
\texttt{url} is a collection of data for idenfitication of malicious URLs.
It has $N = 2396130$ URL samples, each having $d = 64$ real valued features
out of a total of $3231961$ attributes (density: $18.08\%$). Finally,
\texttt{epsilon} is a synthetic, dense dataset, having $N = 500000$ samples and
$d = 2000$ features.

It can be verified that $\nabla F(x)$ is $(1/4 \norm{A}_{2}^{2} +
\lambda_{2})$-Lipschitz continuous with $\norm{A}_{2}^{2} = 1$ in all the examples, and
$F(x)$ is $\lambda_{2}$-strongly convex with respect to $\norm{\cdot}_{2}$.

We create three \texttt{c4.2xlarge} compute nodes in Amazon's Elastic
Compute Cloud. The compute nodes are physically located in Ireland
(\texttt{eu}), North Virginia (\texttt{us}) and Tokyo (\texttt{ap}),
respectively. Then, we assign one CPU from each node as workers, resulting
in a total of 3 workers, and we pick the master node at KTH in Sweden.
We run a small number of iterations of the algorithms to obtain an
\emph{a priori} delay distribution of the workers in this setting, and we
observe that $\bar{\tau} = 6$.

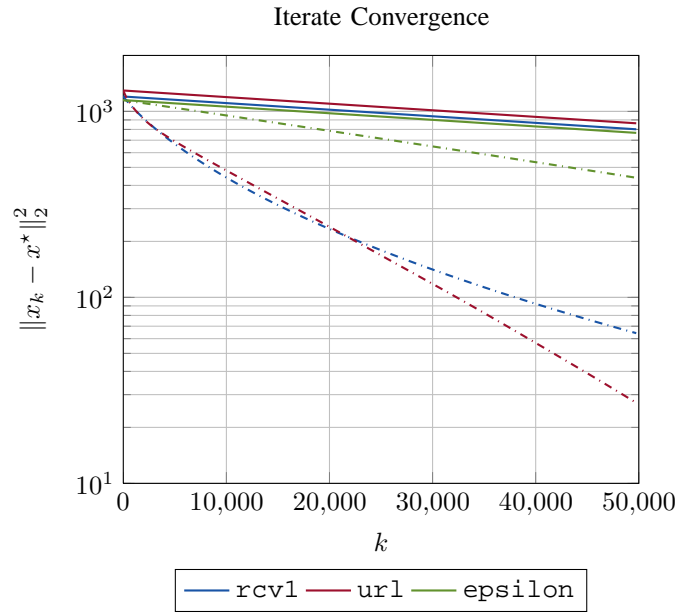
\begin{figure}[ht]
  \centering
  \begin{tikzpicture}
    \begin{semilogyaxis}[%
    xlabel=$k$,%
    ylabel=$\norm{x_{k}-x^{\star}}_{2}^{2}$,%
    grid=both,%
    scaled x ticks=false,%
    xmin=0,xmax=50000,%
    ymin=1e1,ymax=2e3,%
    xtick={0,10000,20000,30000,40000,50000},%
    title=Iterate Convergence,%
    legend style={at={(0.5,-0.20)},anchor=north,legend columns=-1},%
    ]
    % rcv1 bound (blue)
    \addplot[kth-dark-blue,solid,thick]
      table[x index=0,y index=2,col sep=comma,header=true]
      {aws-iterations-rcv1-2-cut.csv};
    % url bound (red)
    \addplot[kth-dark-red,solid,thick]
      table[x index=0,y index=2,col sep=comma,header=true]
      {aws-iterations-url-2-cut.csv};
    % epsilon bound (green)
    \addplot[kth-dark-green,solid,thick]
      table[x index=0,y index=2,col sep=comma,header=true]
      {aws-iterations-epsilon-2-cut.csv};
    % and their experiment results
    \addplot[kth-dark-blue,dashdotted,thick]
      table[x index=0,y index=1,col sep=comma,header=true]
      {aws-iterations-rcv1-2-cut.csv};
    \addplot[kth-dark-red,dashdotted,thick]
      table[x index=0,y index=1,col sep=comma,header=true]
      {aws-iterations-url-2-cut.csv};
    \addplot[kth-dark-green,dashdotted,thick]
      table[x index=0,y index=1,col sep=comma,header=true]
      {aws-iterations-epsilon-2-cut.csv};
    % legend entries
    \legend{\texttt{rcv1},\texttt{url},\texttt{epsilon}}
    \end{semilogyaxis}
  \end{tikzpicture}
  \caption{Convergence of the iterates in Amazon EC2 experiments. Solid lines
    represent our theoretical upper bound, whereas dash-dotted lines represent
    experiment results.}\label{fig:aws-iterates}
\end{figure}

\begin{figure}[ht]
  \centering
  \begin{tikzpicture}
    \begin{axis}[%
    ybar,%
    xlabel=Worker,%
    ylabel=$\tau$,%
    grid=major,%
    cycle list name=kth-colors,%
    symbolic x coords={eu,us,ap},%
    xtick=data,%
    %bar width=8pt,%
    title=Delay Distribution,%
    legend style={at={(0.5,-0.20)},anchor=north,legend columns=-1},%
    ]
    % rcv1 experiments
    \addplot+[error bars/.cd,y dir=both,y explicit]
      table[x index=0,y index=1,y error index=2,col sep=comma,header=true]
      {aws-delays-rcv1-2.csv};
    % url experiments
    \addplot+[error bars/.cd,y dir=both,y explicit]
      table[x index=0,y index=1,y error index=2,col sep=comma,header=true]
      {aws-delays-url-2.csv};
    % epsilon experiments
    \addplot+[error bars/.cd,y dir=both,y explicit]
      table[x index=0,y index=1,y error index=2,col sep=comma,header=true]
      {aws-delays-epsilon-2.csv};
    % legend entries
    \legend{\texttt{rcv1},\texttt{url},\texttt{epsilon}}
    \end{axis}
  \end{tikzpicture}
  \caption{Worker delays in Amazon EC2 experiments. Bars represent the mean
    delays, whereas vertical stacked lines represent the standard deviation.
    For each worker, from left to right, we present the delays obtained in
    \texttt{rcv1}, \texttt{url} and \texttt{epsilon} experiments, respectively.}
    \label{fig:aws-delays}
\end{figure}
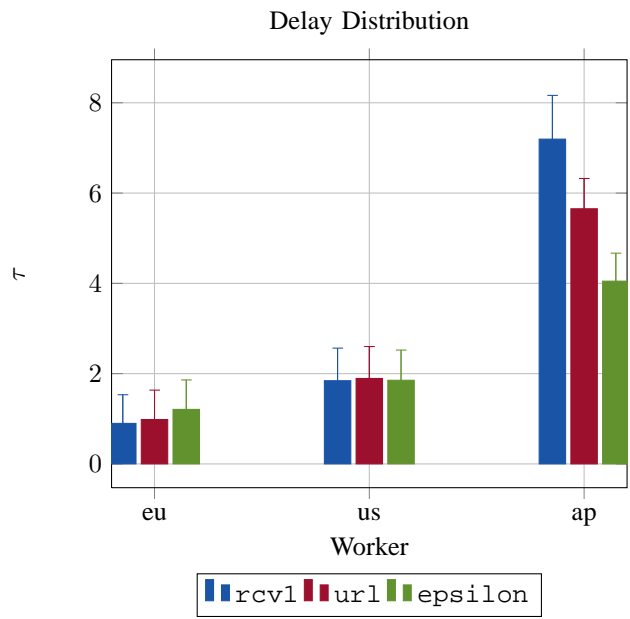

In Figures~\ref{fig:aws-iterates}~and~\ref{fig:aws-delays}, we present the
convergence results of our experiments and delay distributions of the workers,
respectively. As in the previous example, the iterates converge to the optimizer
and the theoretical bound derived is valid. Another observation worth noting is
that the denser the datasets become, the smaller the gap between the actual
iterates and the theoretical upper bound gets.

\section{Discussions and conclusion}

In this paper, we have studied the use of parameter server framework on solving
regularized machine learning problems. One class of methods applicable for this
framework is the proximal incremental aggregated gradient method. We have shown
that when the objective function is strongly convex, the iterates generated by
the method converges linearly to the global optimum. We have also given constant
step-size rule when the degree of asynchrony in the architecture is known.
Moreover, we have validated our theoretical bound by simulating the parameter
server architecture on two different problems.

\bibliographystyle{IEEEtran}
\bibliography{IEEEabrv,references}

\end{document}